\documentclass[12pt,emtex]{report}
\usepackage{times}

\usepackage{graphicx}
\usepackage{amsfonts}       %
\usepackage{amsmath} 
\usepackage{amssymb}        %
\usepackage{amsthm}         

\usepackage{epsfig}
\psfigdriver{dvips}

\usepackage{latexsym}
\usepackage[matrix,curve]{xypic}

\usepackage[german]{babel}  
\newcommand{\A}{\mathbb{A}}

\newcommand{\Q}{\mathbb{Q}}
\newcommand{\R}{\mathbb{R}}

\newcommand{\I}{\mathbb {I}}
\DeclareMathSizes{12}{11}{8}{5}

\begin{document}

\title{The Norm Index Theorem \break (An Analytic Proof) }
\author{Rainer Weissauer}
\maketitle

\bigskip\noindent
\centerline{\bf Introduction}

\bigskip
A key result of class field theory for abelian field extensions $L/K$, concerning the idele
norm map $N:\I_L\to \I_K$, is the norm index theorem
$$[\I_K:K^*N(\I_L)] = [L:K]\ .$$
The inequality $[\I_K:K^*N(\I_L)] \leq [L:K]$ is known to hold in general. A rather easy and
well known analytic proof ([H1]) employs analytic properties of $L$-series near the border
$s=1$ of convergency, obtained by expressing Hecke $L$-series in the form of strongly
converging integrals using the Poisson formula (see [H2], [T1]).

\bigskip
The second inequality $[\I_K:K^*N(\I_L)] \geq [L:K]$ holds for abelian extensions $L/K$ only.
It can be easily reduced to the case of cyclic extensions, for which it is usually proved by
nonanalytic methods (see [T2]).

\bigskip
This note provides a purely analytic proof of the second inequality via the trace formula for
the compact multiplicative group $ L^* Z_L\backslash \I_L$. The Poisson formula, used for the
first inequality, is the trace formula for the  additive compact group $\A_K/K$ enhanced by the
information, that the Pontryagin dual $(\A_K/K)^D$  has dimension 1 as a $K$-vector space.
Hence the spectral theory of the multiplicative and the additive theory combined prove the norm
index theorem by analytic methods.

\bigskip\noindent
\centerline{\bf Review of the trace formula}

\bigskip\noindent
For the number field $L$  let $\I_L$ be the group of ideles. Let denote $Z_L\cong \R^*_{>0}$
the image of $\R^*_{>0}\subseteq \I_\Q$ in $\I_L$. The quotient $X_L=L^*Z_L\backslash \I_L$ is
a compact group. For Haar measures $dg_L$ on $\I_L$ and $dz_L$ on  $Z_L$ the measure
$\frac{dg_L}{dz_L}$ induces a measure $dx_L$ on $X_L$. Assume $\int_{X_L} dx_L=1$. The
corresponding Hilbert space $L^2(X_L)$ is spanned by the characters $\eta\in (X_L)^D$ of $X_L$.
Let $L/K$ be cyclic with Galois group $\langle\sigma\rangle$.

\bigskip
For functions $\prod_v f_v(x_v)$   in $C_c^\infty(\I_L)$ define $f$ by integration over $Z_K$.
For $\varphi\in L^2(X_L)$ define the convolution $R\varphi$ in $L^2(X_L)$ by
$\int_{Z_L\backslash\I_L} f(h^{-1}\theta(g))\varphi(g) \frac{dg_L}{dz_L}\ $ as a function of
$h$. Here $\theta(g)=\kappa \cdot \sigma(g)$ for a fixed $\kappa\in \I_L$ with idele norm
$N(\kappa)\in K^*$. Obviously $\theta=\theta_{\kappa}$ defines an automorphism of $X_L$ of
order $[L:K]$.

\bigskip
The operator $R$ has the kernel $K(y,x)=\sum_{\delta\in L^*} f(y^{-1}\delta \theta(x))$. Hence
its trace is $\int_{X_L} K(x,x) dx_L$. Using $K^*\backslash L^* \cong (\sigma-1)L^*$ and
$Z_L=Z_K$ the integral defining the trace, for $c=\frac{dz_K}{dz_L}$, therefore becomes
$$ \sum_{\delta\in L/(\sigma -1)L^*} \int_{y\in \I_K\backslash \I_L}  \Bigl( c \cdot \int_{x\in Z_K
K^*\backslash\I_K} dx_K \Bigr) f\Bigl(y^{-1}\delta \theta(y)\Bigr)\ \frac{dg_L}{dg_K}  \ .$$
Abbreviate $O^L_\delta(f_w)= \int_{K_w^*\backslash L_v^*}\ \bigl(\prod_{v\vert w}
f_v(g_v^{-1}\delta \theta(g_v)dg_v\bigr)\ /dg_w$, so this simplifies to
$$ \sum_{\delta\in L^*/(\sigma -1)L^*} c \cdot \prod_v \ O^L_\delta(f_v)\ .$$ For characters $\eta$ of
$X_L$ put $\eta^\theta(x)=\eta(\sigma^{-1}(x))$ and $\eta(f)=\int_{Z_L\backslash\I_L}
f\bigl(\theta(g)\bigr)\eta(g)\frac{dg_L}{dz_L}$. Hence up to these constants $\eta(f)$ $$
R\eta(h)\ = \ \eta(f)\cdot \eta^\theta(h)\ .$$ In other words, the trace of $R$ becomes  the
spectral sum $ \sum_{\eta=\eta^\theta} \eta(f) $. Comparing with the previous formula we obtain
the usual {trace formula} (as in [KS])
  $$\quad \sum_{\eta=\eta^\theta}\ \eta(f)\ = \sum_{\delta\in
L^*/(\sigma -1)L^*} \ c \cdot \prod_w \ O^L_\delta(f_w) \ .$$

\bigskip\noindent
\centerline{\bf Matching functions}

\bigskip\noindent
$\I_K =\prod_w K_w^*$ and $\I_L=\prod_w L_w^*=\prod_w\prod_{v\vert w} L_v^*$. Functions
$\prod_w f_w(x)$ in $C_c^\infty(\I_L)$, s.t. $f_w=\prod_{v\vert w} f_v$, and functions $\prod_w
h_w(x)$ in $C_c^\infty(\I_K)$ are said to be matching functions, if
$h_w(\gamma_w)=O^L_{\delta_w}(f_w)$ holds for $\gamma_w=N(\delta_w),\delta_w\in L_w^*$ and
$h_w(\gamma_w)$ is zero for $\gamma_w\notin N(L_w^*)$. Notice, that $h_w$ is uniquely
determined by $f_w$. Existence is obvious, since the characteristic functions $\prod_{v\vert
w}1_{{\frak o}_v}$ and $1_{{\frak o}_w}$ of integral elements do match at all unramified
nonarchimedean places by the elementary property $N({\frak o}_v^*)={\frak o}_w^*$, which is
valid for all unramified  places $w$ (the fundamental lemma).

\goodbreak
\bigskip\noindent
\centerline{\bf Twisted case revisited}

\bigskip\noindent
Characters $\eta\in (X_L)^D$ on the spectral side of the trace formula are characters
$\eta=\eta^\theta$ of $X_L$ trivial on $(\sigma -1)X_L$, hence
 of the form $\eta=
\chi^\sharp \circ N$ for characters $\chi^\sharp$ of $Y^\sharp= N(Z_L\backslash
\I_L)\big/N(L^*)=X_L/(\sigma -1)X_L$ (Hilbert theorem 90). For $\kappa=1$ the trace summands
$\eta(f)$ therefore can be written in the form
$$ \eta(f)= \int_{Z_LKern(N)\backslash \I_L} \eta(g)\ \Bigl(\int_{Kern(N)} f(gn)dn\Bigr)\ \frac{dg_L}{dn dz_L}
\ $$ for $Kern(N)=(\sigma -1)\I_L$ and $dn= \frac{dg_L}{dg_K}(h)$ and $n=\sigma(h)h^{-1}$. By
the matching condition and $\eta(g)=\chi^\sharp\bigl(N(g)\bigr)$ the last expression giving
$\eta(f)$ becomes
$$ \tilde c \cdot \chi^\sharp(h)= \int_{Z_K\backslash N(\I_L)} \tilde c \cdot \chi^\sharp\bigl(N(g)\bigr)\ h\bigl(N(g)\bigr)\ \frac{dg_K}{dz_K} \quad \mbox{ where } \quad  \tilde c\cdot \frac{dg_K}{dz_K} =\frac{dg_L}{dn dz_L}
\ .$$ On the right side of the  trace formula we can also apply the matching condition. Hilbert
90 implies $L^*/(\sigma -1)L^*\cong N(L^*)$, hence
 the still preliminary formula $(\sharp)$
$$ \sum_{\chi^\sharp \in (Y^\sharp)^D}\ \tilde c \cdot \chi^\sharp(h)\ = \sum_{\gamma\in N(L^*)}\ c\cdot h(\gamma) \ .$$

\goodbreak
\bigskip\noindent
\centerline {\bf The degenerate case $K=L$}

\bigskip\noindent
For a pair of matching functions we compare the last  formula $(\sharp)$ with the trace formula
for $h$ in the case $L=K$ and $\kappa=1$, which simply reduces to
$$ \sum_{\chi\in (X_K)^D}\ \chi(h) = \sum_{\gamma\in K^*} h(\gamma) \ ,$$
since $\prod_w O^K_\gamma(h_w)=h(\gamma)$. By the matching condition the support of  $h$ is
contained in $N(Z_L\backslash \I_L)\subseteq Z_K\backslash \I_L$. Therefore on the left we may
restrict characters $\chi$ from $X_K$ to the image $Y^\flat=N(Z_L\backslash \I_L)\big/(K^*\cap
N(\I_L))$ of $N(Z_L\backslash \I_L)$ in $X_K $, which is a subgroup of index $[\I_K: K^*
N(\I_L)]$. For the restrictions  $\chi^\flat$ of the characters $\chi$, now with $\chi^\flat$
running over the character group $(Y^\flat)^D$ of $Y^\flat$, we thus may restate the trace
formula for $L=K$ as the following formula $(\flat)$
$$ \#\Bigl(\frac{\I_K}{ K^* N(\I_L)}\Bigr)\cdot \sum_{\chi^\flat\in (Y^\flat)^D}\ \chi^\flat(h)\ = \sum_{\gamma\in K^*\cap N(\I_L)} h(\gamma) \ .$$

In particular $[\I_K: K^* N(\I_L)]<\infty$.

\goodbreak
\bigskip\noindent
\centerline{\bf Comparing the trace formulas}

\bigskip\noindent
Recall $Y^\sharp= X_L/(\sigma -1)X_L =N(Z_L\backslash\I_L)\Big/N(L^*)$, which  gives the exact
sequence
$$ 0 \to \frac{K^*\cap N(\I_L)}{N(L^*)} \to Y^\sharp \to Y^\flat \to 0 \ .$$
A system of representatives $\kappa_i$ for all possible  $\kappa$ modulo $L^*$ is in 1-1
correspondence with the elements of $\frac{K^*\cap N(\I_L)}{N(L^*)}$. Summing up the
$\kappa_i$-twisted trace formulas for $f$ -- these are nothing but the revisited forms of the
trace formulas $(\sharp)$ for the translates $f(\kappa_i x)$ of $f(x)$ -- we obtain from
$(\sharp)$ therefore the final identity
$$ \#\Bigl(\frac{K^*\cap N(\I_L)}{N(L^*)}\Bigr) \cdot \sum_{\chi^\flat \in (Y^\flat)^D}\ \tilde c \cdot  \chi^\flat(h)\ =  \sum_{\gamma\in K^*\cap
N(\I_L)} c\cdot h(\gamma) \ .$$

In particular $[K^*\cap N(\I_L) : N(L^*)] <\infty$. If we compare this last formula with the
trace formula $(\flat)$ for $L=K$, we get the crucial formula
$$ \frac{[\I_K : K^* N(\I_L)]}{[K^*\cap N(\I_L) : N(L^*)]}\ =\ \frac{\tilde c}{c}\ .$$

\bigskip
{\it The quotient $\tilde c/c$}.
The constants were defined by $\tilde c \cdot \frac{dg_K}{dz_K} = \frac{dg_L}{dn dz_L}$ and $c=
\frac{dz_K}{dz_L}$, and $dn= \frac{dg_L}{dg_K}$ by abuse of notation.  Indeed, the ratio
$\tilde c/c$ is independent from the particular choice of Haar measures $dg_L,dg_K$ and
$dz_K,dz_K$. Therefore we may choose them freely. Normalizing constants for $dg_K,dg_L$ cancel.
Unraveling the definitions in terms of the maps $i$ and $N$ we are thus reduced to consider
invariant $K$-rational differential forms for the $K$-tori $\mathbb G_m$ and
$T=Res_{L/K}({\mathbb G}_m)$ and the exact sequence
$$ 0 \to {\mathbb G}_m \overset{i}{\to} T \overset{1-\sigma}{\to} T \overset{N}{\to} {\mathbb G}_m \to 0 \ .$$
An easy calculation on the tangent spaces, for $A=i^*(e_1^*)$ and a form $B$ on the tangent
space of $V=(1-\sigma)T$, using the formula  $e_1^*\wedge (1-\sigma)^*(B) = e_1^*\wedge
[(e_2^*-e_1^*)\wedge \cdots (e_n^*-e_{n-1}^*)] = e_1^*\wedge \cdots \wedge e_n^* = [e_1^*\wedge
\cdots e_{n-1}^*] \wedge (e_1^*+\cdots + e_n^*)= B\wedge N^*(A)$ proves, that the quotient
$\tilde c/c$ entirely comes from the measure comparison between $dz_L$ and $dz_K$ similarly
arising from the exact sequence
$$ 0 \to Z_K \overset{i}{\to} Z_L \overset{1-\sigma}{\to} Z_L \overset{N}{\to} Z_K \to 0 \ .$$
$(1-\sigma)$ is the zero map on $Z_L$. Hence this comparison is trivial. The factor $\tilde
c/c$ immediately turns out to be $\tilde c/c=[L:K]$. This completes the proof. Of course,
combined with the first inequality, this a posteriori implies the Hasse norm theorem $K^*\cap
N(\I_L) = N(L^*)$, hence the stability of our particular trace formula as in [KS] 6.4 and 7.4.

\bigskip\noindent
\centerline{\bf Bibliography}

\bigskip
{\sc Hecke E.}, [H1] {\it \"Uber eine neue Anwendung der Zetafunktionen auf die Arithmetik der
Zahlk\"orper} (1917), Gesammelte Werke n.8 (1959), page 172 -177, and [H2] {\it \"Uber die
Zetafunktion beliebiger algebraischer Zahlk\"orper} (1917), Gesammelte Werke n.7 (1959), page
159 -171

\bigskip
{\sc Kottwitz R.E.-Shelstad D.}, [KS] {\it Foundations of Twisted Endoscopy}, asterisque 255,
(1999)

\bigskip
{\sc Tate J.T.}, [T1] {\it Fourier Analysis in Number Fields and Hecke's Zeta-Functions}, and
[T2] {\it Global Class Field Theory}, in {\it Algebraic Number Theory}, ed. J.W.S. Cassels and
A. Fr\"ohlich (1967), Academic Press

\end{document}